\documentclass[12pt,a4paper]{elsarticle}
\usepackage[T2A]{fontenc}
\usepackage{amssymb, latexsym, amsmath}
\usepackage[english]{babel}
\usepackage{graphicx}
\usepackage{amsfonts,amsmath,color,verbatim}
\usepackage[framemethod=tikz]{mdframed}
\usepackage{lipsum}
\usepackage{float} 
\usepackage{subcaption}
\usepackage{tikz}
\usetikzlibrary{patterns}
\usetikzlibrary{arrows,decorations,backgrounds}
\usetikzlibrary{decorations.markings}
\usetikzlibrary{arrows.meta}

\newcommand{\btau} {\mbox{\boldmath $\tau$}}
\newcommand{\bgam} {\mbox{\boldmath $\gamma$}}

\begin{document}
\begin{frontmatter}

\title{ Squeeze flow of a Hershel--Bulkley fluid}

\author{Larisa Muravleva}

\address{Lomonosov Moscow State University}
 
\ead{lvmurav@gmail.com}

\begin{abstract}
We develop an asymptotic solution for the axisymmetric squeeze flow of a viscoplastic Hershel--Bulkley medium following the asymptotic technique suggested earlier by Balmforth and Craster (1999) and Frigaard and Ryan (2004).
\end{abstract}

\end{frontmatter}

\section{Introduction} 

Viscoplastic or yield-stress fluids are materials which behave like a solid below  a critical
yield stress and flow like a viscous fluid for stresses higher than this threshold.
The flow field is thus divided
into unyielded (rigid) and yielded (fluid) zones. The surface separating a rigid from a fluid zone is known as a yield surface. The location and shape of the latter must be determined as part of the solution of the flow problem. 

The squeeze flow of viscoplastic materials has been investigated many times, see the detailed reviews \cite{Eng}-- \cite{MitTsam}
for a comprehensive list of references to available analytical, numerical and experimental results \cite{CovStan}--\cite{Lawal98}
and \cite{AlGeor2013}--\cite{FusFar} for some
more recent contributions.
Different constitutive equations have been used, in both theoretical and
numerical studies, such as the original Bingham model by Covey and Stanmor \cite{CovStan}, Lipscomb and
Denn\cite{LipDenn}, Sherwood and Durban \cite{SherDur96}, Smyrnaios and
Tsamopoulos \cite{ST}, and Roussel et al. \cite{Lanos06}, the bi-viscosity model by
 O'Donovan and Tanner \cite{DonTan}] and Wilson \cite{Wil93}, the Hershel-Bulkley 
model by Covey and Stanmor \cite{CovStan}, Sherwood and Durban \cite{SherDur98}, an elasto-viscoplastic model by Adams et al. \cite{Adams97}, Bingham fluid with a deformable core by Fusi et al. \cite{FusFar} and the regularized Papanastasiou
 model by Smyrnaios and Tsamopoulos \cite{ST}, Matsoukas
and Mitsoulis \cite{MatMit}, \cite{MitMat2005}, and Karapetsas
and Tsamopoulos \cite{KarTsam}.

 We would like to construct the consistent asymptotic solution for the squeeze flow of Hershel-Bulkley material.
  Walton and Bittleston \cite{WaltBit} studied the axial flow in an eccentric annular duct and showed analytically that a true plug exists in the middle of the channel on both the wide and narrow part of the annulus, with a pseudo-plug placed between the two rigid zones. Balmforth and Craster \cite{BalmCr}, Frigaard and Ryan \cite{FrigRyan}, \cite{PUTZ}
suggested the asymptotic technique that allows constructing the consistent solution for thin--layer problems. The asymptotic solution was developed for a fluid flowing down an inclined plane \cite{BalmCr} and for the flow along a channel of slowly varying width \cite{FrigRyan}, \cite{PUTZ}. 
Recently, Muravleva \cite{LMur}, \cite{Mur_JNNFM2017} has analysed the planar and axisymmetric squeeze flows of a Bingham fluid and the axisymmetric squeeze flows of a Casson fluid \cite{Mur_JNNFM2017} exploiting the asymptotic technique introduced in \cite{BalmCr}--\cite{PUTZ}.
In this article, we continue our research into the squeeze flow problem \cite{LMur}-\cite{Mur_JNNFM2017} and consider the axisymmetric squeeze flow of a Hershel-Bulkley material, following an approach developed in \cite{BalmCr}--\cite{PUTZ}. 

The paper is organized as follows. In Section 2 the dimensionless
governing equations of the flow are presented. Section 3 is the 
Finally, in Section 5, a summary of the results is given.

\section{Problem statement.}

\begin{figure}
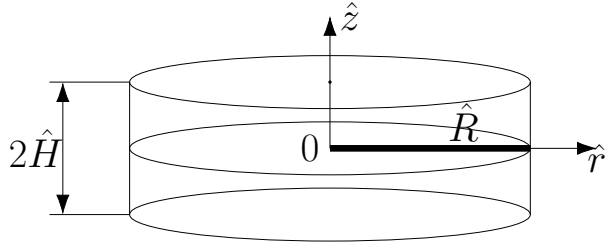

\centering
\resizebox{.6\textwidth}{!}{
\tikz{
\draw (0,2) ellipse (3cm and .4cm);
\draw (0,1) ellipse (3cm and .4cm);
\draw (0,0) ellipse (3cm and .4cm);
\draw (-3,0) -- (-3,2);
\draw (3,0) -- (3,2);
\node [xshift=-0.3cm,yshift=1.0cm] {\Large 0};
\fill (0,2) circle (0.025);
\draw [-{Latex[length=3mm, width=2mm]}](0,1) -- (4,1);
\node [xshift=4cm,yshift=0.85cm] {\Large $\hat{r}$};
\draw[-{Latex[length=3mm, width=2mm]}] (0,1)--(0,3);
\node [xshift=0.3cm,yshift=3cm] {\Large $\hat{z}$};
\draw[{Latex[length=3mm, width=2mm]}-{Latex[length=3mm, width=2mm]}] (-4,0)--(-4,2);
\node [xshift=-4.4cm,yshift=1cm] {\Large $2 \hat{H}$};
\draw (-3,0) -- (-4.2,0);
\draw (-3,2) -- (-4.2,2);
\draw [line width=0.1cm](0,1) -- (3,1);
\node [xshift=2cm,yshift=1.4cm] {\Large $\hat{R}$};
}}
\caption{Coordinate system and basic dimensions used to describe axisymmetric squeeze flows.}
\label{mur_fig1}
\end{figure}

Squeeze Flow (SF) is the process in which a fluid is squeezed between two
approaching parallel plates resulting in a radial flow, outward from the center.
The geometry of the problem is shown in Fig.~\ref{mur_fig1}: we use an axisymmetric
cylindrical polar coordinate system $(r; \theta; z)$ to describe the squeezing of a
cylinder of incompressible Hershel-Bulkley fluid with radius $\hat{R}$ and height $2 \hat{H}$. The
fluid has density $\hat{\rho}$, yield stress $\hat{\tau_0}$ and plastic viscosity $\hat{\mu}$. The plates are
squeezed together at a velocity $\hat{W}$, inducing a flow. We denote the dimensional
variables with a hat symbol.
We have scaled lengths in the $r$ and $z$ directions differently,
with the disk radius $\hat{R}$ and
with the half-distance $\hat{H}$ between the disks, respectively. $\hat{W}$ is taken as the characteristic velocity in the transverse direction, and the radial velocity
component is scaled with $\hat{U} = \hat{W}\hat{R}/\hat{H}$.
 The pressure is scaled with $\hat{\mu} \hat{W}^n \hat{R}^{n+1}/ \hat{H}^{2n+1}$, and time with $\hat{H} / \hat{W}$. The shear-stress components are scaled with $\hat{\mu} \hat{W}^n \hat{R}^n/\hat{H}^{2n}$, the extensional stresses with $\hat{\mu} \hat{W}^n \hat{R}^{n-1} / \hat{H}^{2n-1}$. 
 
The flow is governed by the dimensionless conservation equations of
momentum and mass:
\begin{equation}\label{mur_1}
\varepsilon Re \Bigl(\frac{\partial u}{\partial t} +
u \frac{\partial u}{\partial r} + w \frac{\partial u}{\partial z}\Bigr) =
- \frac{\partial p}{\partial r} + \varepsilon^2 \frac{\partial \tau_{rr}}{\partial r} +
\frac{\partial \tau_{rz}}{\partial z} + \varepsilon^2 \frac{\tau_{rr} - \tau_{\theta\theta}}{r},
\end{equation}
\begin{equation}\label{mur_2}
\varepsilon^2 Re \Bigl( \frac{\partial w}{\partial t} +
u \frac{\partial w}{\partial r} + w \frac{\partial w}{\partial z} \Bigr) =
- \frac{\partial p}{\partial z} + \varepsilon^2 \frac{\partial \tau_{rz}}{\partial r} - \varepsilon^2
\frac{\partial \tau_{rr}}{\partial z}- \varepsilon^2 \frac{\partial \tau_{\theta\theta}}{\partial z}
+\varepsilon^2 \frac{\tau_{rz}}{r},
\end{equation}
\begin{equation}\label{mur_3}
\frac{\partial u}{\partial r} + \frac{u}{r} + \frac{\partial w}{\partial z} = 0.
\end{equation}
The Herschel-Bulkley plastic stresses are related to the strain rates through
the constitutive equations
\begin{eqnarray}\label{mur_4}
\sigma_{ij} = -p \delta_{ij} + \tau_{ij},\quad
\begin{cases}
\tau_{ij} = ( \dot{\gamma}^{n-1} + \frac{B}{\dot{\gamma}})\dot{\gamma}_{ij},   &\mbox{ if }\tau > B,\\
\dot{\gamma}_{ij} = 0, &\mbox{ if }\tau \le B.
\end{cases}
\end{eqnarray}
where ${\tau}=\sqrt{\frac{1}{2}\btau:\btau}$ and $\dot{\gamma}=\sqrt{\frac{1}{2}\dot{\bgam}:\dot{\bgam}}$ denote the second invariants of
${\btau}$ and ${\dot{\bgam}}$, and the strain rate tensor ${\dot{\bgam}}$ is 
given by
\begin{equation}\label{mur_6}
{\dot{\gamma}}_{rr} = 2 \frac{\partial u}{\partial r},\quad
{\dot{\gamma}}_{rz} = \Bigl( \frac{\partial u}{\partial z} + \varepsilon^2 \frac{\partial w}{\partial r}\Bigr),\quad
{\dot{\gamma}}_{zz} = 2 \frac{\partial w}{\partial z},\quad
{\dot{\gamma}}_{\theta\theta} = 2 \frac{u}{r},
\end{equation}
The scalings introduce
the small aspect ratio 
$\varepsilon$, and the
Reynolds and Bingham numbers:
\begin{equation}\label{mur_5}
\varepsilon = \frac{\hat{H}}{\hat{R}},\quad
Re = \frac{\hat{\rho}\hat{W}\hat{R}}{\hat{\mu}},\quad B = \frac{\hat{\tau_0}\hat{ H}}{\hat{\mu}\hat{U}}.
\end{equation}

We consider $Re \ll 1 $. Neglecting the inertial terms, \eqref{mur_1}--\eqref{mur_3} are replaced by:
\begin{equation}\label{mur_10}
- \frac{\partial p}{\partial r} + \varepsilon^2 \frac{\partial \tau_{rr}}{\partial r} +
\frac{\partial \tau_{rz}}{\partial z} + \varepsilon^2 \frac{\tau_{rr} - \tau_{\theta\theta}}{r} = 0,
\end{equation}
\begin{equation}\label{mur_11}
- \frac{\partial p}{\partial z} + \varepsilon^2 \frac{\partial \tau_{rz}}{\partial r} - \varepsilon^2
\frac{\partial \tau_{rr}}{\partial z}- \varepsilon^2 \frac{\partial \tau_{\theta\theta}}{\partial z}+
\varepsilon^2 \frac{\tau_{rz}}{r} = 0,
\end{equation}
\begin{equation}\label{mur_12}
\frac{\partial u}{\partial r} + \frac{u}{r} + \frac{\partial w}{\partial z} = 0.
\end{equation}
Exploiting the symmetry of the axisymmetric flow 
about the plane $z=0$,
we solve these equations over the domain $0\le z \le 1$, $0\le r \le 1$ subject to no-slip conditions, $u = u_b$, $w = -1$ on the surface of the disc $z = 1$, symmetry conditions $\tau_{rz} = 0$, $w = 0$ on the plane symmetry $z = 0$ and $u = 0$, $\tau_{rz} = 0$ on the axis symmetry $r = 0$, and stress-free 
$\sigma_{rr} = - p + \varepsilon^2 \tau_{rr} = 0$, $\tau_{rz} = 0$ at the edge $r = 1$. 

\section{Asymptotic expansions.}\label{sec-asymp}

We now solve the equations by introducing an asymptotic expansion.
First, we consider shear flow near the plate
for which we may find a solution through
a straightforward expansion of the equations. This shear solution is denoted
with a superscript (s).

We consider regular expansions in $\varepsilon$ of form:
\begin{equation}\label{ser_sh1}
\begin{split}
&u^s = u^{s,0} + \varepsilon u^{s,1} + \varepsilon^{s,2} u^2 \ldots,\quad
w^s = w^{s,0} + \varepsilon w^{s,1} + \varepsilon^2 w^{s,2} \ldots,\quad\\
&p^s = p^{s,0} + \varepsilon p^{s,1} + \varepsilon^2 p^{s,2} \ldots,\quad\tau^s_{ij} = \tau_{ij}^{s,0} + \varepsilon \tau_{ij}^{s,1} + \varepsilon^2 \tau_{ij}^{s,2} \ldots.
\end{split}
\end{equation}
We substitute these expansions into the governing equations \eqref{mur_10}--\eqref{mur_12}, 
and collect together the terms of the same order.

 The lubrication equations for the first two orders are:
\begin{eqnarray}
\label{sh1}
&\mathcal{O}(1)\quad\quad
&- \frac{\partial p^{s,0}}{\partial r} + \frac{\partial \tau_{rz}^{s,0}}{\partial z} = 0,\\
\label{sh2}
&&- \frac{\partial p^{s,0}}{\partial z} = 0,\\
\label{sh3}
&&\frac{\partial u^{s,0}}{\partial r} + \frac{u^{s,0}}{r} + \frac{\partial w^{s,0}}{\partial z} = 0.
\end{eqnarray}
\begin{eqnarray}
\label{first_1}
&\mathcal{O}( \varepsilon)\quad\quad
&-\frac{\partial p^{s,1}}{\partial r} + \frac{\partial \tau_{rz}^{s,1}}{\partial z} = 0,\\
\label{first_2}
&&-\frac{\partial p^{s,1}}{\partial z} = 0,\\
\label{first_3}
&&\frac{u^{s,1}}{r} +\frac{\partial u^{s,1}}{\partial r} + \frac{\partial w^{s,1}}{\partial z} = 0.
\end{eqnarray}
From the
conservation of mass we have $\int_{0}^{1} u (r, z)~dz = \frac{r}{2}$. We require that
\begin{equation}\label{q}
\int_{0}^{1} u^0 (r, z)~dz = \frac{r}{2}, \quad \int_{0}^{1} u^1 (r, z)~dz = 0.
\end{equation}

\subsection{Zero-order approximation}\label{sec-shear}
After the solution of the equations  \eqref{sh1}, \eqref{sh2} we have
\begin{equation}\label{p_0}
p^{s,0} = p_0 (r),\quad \tau_{rz}^{s,0} = zp_0'(r),
\end{equation}
where  $p_0$ is a function only of $r$ and the prime $( )'$ represents  derivative ${\frac{d}{dr}}$.
The leading order second invariants of strain rate and stress
are given by
$\tau^{s,0} = |\tau^{s,0}_{rz}|$,
$\dot{\gamma}^{s,0} = \Bigl| {\frac{\partial u^{s,0}}{\partial z}} \Bigr|$.
We are looking for a solution with $\frac{\partial u^{s,0}}{\partial z} < 0$ in the domain $r >0$, $z > 0$.
Provided the yield stress is exceeded, the stress tensor components become
\begin{equation}\label{sh5}
\tau_{rz}^{s,0} =  -\Bigl(\Bigl|\frac{\partial u^{s,0}}{\partial z}\Bigr|^n+ B\Bigr), \quad
\tau_{rz}^{s,1} = n \Bigl|\frac{\partial u^{s,0}}{\partial z}\Bigr|\frac{\partial u^{s,1}}{\partial z}.
\end{equation}

We note that the leading order second invariant of stress $\tau^{s,0}
= |\tau_{rz}^0| = z |p_{0}'(r)|$, so
$\tau^{s,0}$ attains its maximum at $z = 1$ and vanish at $z = 0$.
Therefore there exists the point $z = z_0$ at which $\tau^{s,0} = B$ and
$\dot{\gamma}^0 = 0$. Hence, at leading order, the yield condition holds at this point and $z_0 (r) = \frac{B}{|p_0'(r)|}$
is the position of the pseudo-yield surface. For  $z\in [0,\,z_0]$ we have $\tau^{s,0} < B$ and $\dot{\gamma}^{s,0}=\frac{\partial u^{s,0}}{\partial z} = 0$. The velocity field is now given by
\begin{equation}\label{u_0}
u^0(r,z) =
\begin{cases}
\Bigl(\frac{B}{z_0}\Bigr)^{1/n}\frac{1}{1+1/n}[(1-z_0)^{1+1/n}-
(z-z_0)^{1+1/n}], & z \in [0,\,z_0],\\
\Bigl(\frac{B}{z_0}\Bigr)^{1/n}\frac{1}{1+1/n}(1-z_0)^{1+1/n}, & z \in (z_0,\,1].
\end{cases}
\end{equation}
We denote the pseudo-plug velocity by $u_0(r)$: 
\begin{equation}\label{def_u0}
u_0(r) = \frac{B^{1/n}(1-z_0)^{1+1/n}}{z_0^{1/n}(1+1/n)}.
\end{equation}
The expressions \eqref{p_0} can be written as
\begin{equation}\label{mur_tau_s}
 p^{s,0} = p_0 (r),\quad p_0 ' = - \frac{B}{z_0(r)},\quad \tau_{rz}^{s,0} = - \frac{Bz}{z_0(r)}.
\end{equation}

Substituting \eqref{u_0} into the first equation of the flow rate constraint \eqref{q} leads to the following equation 
 where the unknown is the pseudo-yield surface $z_0 (r)$
\begin{equation}\label{buck}
\frac{(1-z_0)^{2+1/n}}{(2+1/n)}-
(1-z_0)^{1+1/n} + \Bigl(\frac{z_0}{B}\Bigr)^{1/n}\frac{r}{2}(1+1/n) = 0.
\end{equation}

\begin{figure}[htb]%
\includegraphics[height=0.4\hsize,width=0.8\hsize]{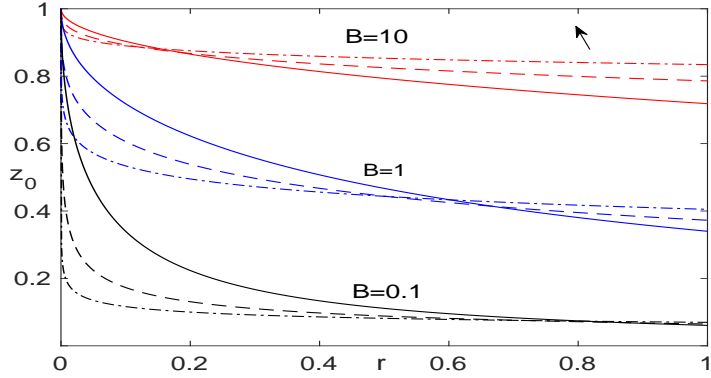}
\caption{
The pseudo-yield surface ${z_0}(r)$ \eqref{buck} for $n=1$ (solid lines), $n=0.5$ (dashed lines), $n=0.25$ (dash-dotted lines).} 
\label{mur_fig2}
\end{figure}

Evidently, $z_0(r)$ and therefore $u_0(r)$ depend on position $r$, highlighting how the flow
in $0 < z \le z_0$ is only a pseudo-plug, which is an extension in the  radial direction.
Thus, the pseudo-yield surface $z_0 (r)$ separates the entire area occupied by the material into subregions: fully yielded zones located near the plates (shear region), and pseudo-plug containing the central plane $z = 0$. 
In the pseudo-plug, the asymptotic expansion \eqref{ser_sh1} breaks down. Therefore, we look for a slightly different asymptotic expansion of the radial velocity component, where the property $\frac{\partial u^{p,0}}{\partial z} = 0$ at $z_0$
is explicitly built in:
\begin{equation}\label{mur_ser_u}
u^p = u^{p,0} (r) + \varepsilon u^{p,1} (r,z) + \varepsilon^2 u^{p,2}(r,z) + \ldots,
\end{equation}
The pseudo-plug solution is denoted with a superscript (p).
The stresses are now given by
\begin{equation}\label{tau}
\tau_{rr}^{p,-1} = \frac{2B}{\dot{\gamma}^{p,0}}\frac{\partial u^{p,0}}{\partial r},\quad
\tau_{rz}^{p,0} = \frac{B}{\dot{\gamma}^{p,0}}\frac{\partial u^{p,1}}{\partial z} ,\quad
\tau_{\theta\theta}^{p,-1} = \frac{2B}{\dot{\gamma}^{p,0}} \frac{u^{p,0}}{r},
\end{equation}
\begin{equation}\label{gam0}
\mbox{where}\quad \dot{\gamma}^{p,0} = \sqrt{
\Bigr( \frac{\partial u^{p,1}}{\partial z}\Bigl)^2
+ 4 \Bigl[\Bigl(\frac{\partial u^{p,0}}{\partial r}\Bigl)^2 +
\Bigl(\frac{u^{p,0}}{r}\Bigr)^2 + \frac{\partial u^{p,0}}{\partial r}\frac{u^{p,0}}{r}\Bigr]
}.
\end{equation}
 The lubrication equations for the first two orders are:
\begin{eqnarray}
\label{pp1}
&\mathcal{O}(1)\quad\quad
&- \frac{\partial p^{p,0}}{\partial r} + \frac{\partial \tau_{rz}^{p,0}}{\partial z} = 0,\\
&&\label{pp2}
- \frac{\partial p^{p,0}}{\partial z} = 0,\\
\label{pp3}
&&\frac{\partial u^{p,0}}{\partial r} + \frac{u^{p,0}}{r} + \frac{\partial w^{p,0}}{\partial z} = 0.
\\
\label{pl1}
&\mathcal{O}( \varepsilon)\quad\quad
&-\frac{\partial p^{p,1} }{\partial r} + \frac{\partial \tau_{rr}^{p,-1}}{\partial r} +
\frac{\tau_{rr}^{p,-1} - \tau_{\theta\theta}^{p,-1}}{r} + \frac{\partial \tau_{rz}^{p,1}}{\partial z} = 0,\\
\label{pl2}
&&-\frac{\partial p^{p,1}}{\partial z} - \frac{\partial (\tau_{rr}^{p,-1} + \tau_{\theta\theta}^{p,-1})}{\partial z} = 0,\\
\label{pl3}
&&\frac{\partial u^{p,1}}{\partial r} + \frac{u^{p,1}}{r}  + \frac{\partial w^{p,1}}{\partial z} = 0.
\end{eqnarray}
Using symmetry about the center plane, we have
$w^{p} \left.\right|_{z = 0} = 0$, $\tau_{rz}^p \left.\right|_{z = 0} = 0$,
 $\frac{\partial u^p}{\partial z} \left.\right|_{z = 0} = 0$.
Integrating  \eqref{pp1}, \eqref{pp2}  and enforcing continuity $p^{0}$,
$\tau_{rz}^0$, $u^0$ at $z = z_0(r)$ leads to
\begin{equation}\label{mur_pl1}
 p^{p,0} (r) = p_0 (r),\quad p_0 ' = - \frac{B}{z_0},\quad \tau_{rz}^{p,0} = - \frac{B z}{z_0},\quad
u^{p,0} = u_0.
\end{equation}
Fig.~\ref{mur_fig2} shows the graphs of the pseudo-yield surfaces $z = z_0(r)$, for different values of $B$ and $n$.
We see that $z = z_0(r)$ is the decreasing function of radius. This property is easy to confirm analytically, differentiating  $z_0$ \eqref{buck} with respect to $r$ ($0 <z_0 \le 1$)
\begin{equation}\label{dz0ns}
z_0'=-\frac{(n+1)(2n+1)z_0^{1+1/n}}{2 B^{1/n}(1-z_0)^{1/n}(2n^2z_0^2 +2n z_0+1+n)}< 0.
\end{equation}
Graphs corresponding to large values of the Bingham number (for fixed values of power-law index $n$) are located above. For a fixed  Bingham number and variable power-law index $n$, the pattern is more complicated: for smaller values of $n$, the graphs are located higher near the center of the plate, and quickly decrease with increasing $r$.
\subsection{The first-order approximation.}\label{sec_first_ns}
In the shear region we integrate \eqref{first_1}, \eqref{first_2}, \eqref{sh5},
using $\tau_{rz}(r,0) = 0$, and receive:
\begin{equation}\label{sh_1ord}
p^{s,1} = p_1 (r),\quad
\tau_{rz}^{s,1} = z p_{1}'(r) + g (r),
\end{equation}
\begin{equation}\label{u_1s}
u^{s,1} = \frac{p_1'(r)}{n+1}\Bigl[(z-z_0)^{\frac{1}{n}}(z+n z_0) - (1-z_0)^{\frac{1}{n}}(1+n z_0)\Bigr]+g(r) \Bigl[(z-z_0)^{\frac{1}{n}} - (1-z_0)^{\frac{1}{n}}\Bigr],
\end{equation}
where  $p_1$ is a function only of $r$ and $g (r)$ is an unknown function of integration.

In the pseudo-plug  region
from \eqref{tau} and \eqref{mur_pl1} we obtain the second invariant $\tau^{p,-1}$ of the stress tensor, which is equal to $B$ since the pseudo-plug region is just at the point of yielding:
\begin{eqnarray}\label{tau0}
&&(\tau^{p,-1})^2 =(\tau_{rr}^{p,-1})^2 + (\tau_{\theta\theta}^{p,-1})^2 +
\tau_{rr}^{p,-1} \tau_{\theta\theta}^{p,-1} + (\tau^{p,0}_{rz})^2 
\notag\\
&& = \frac{4 B^2}{( \dot{\gamma}^{p,0} )^2}
 \Bigl[ (u_0')^2 + \bigl(\frac{u_0}{r}\bigr)^2 + u_0'(r) \frac{u_0}{r} \Bigr] +
 \frac{B^2 z^2}{z_0^2} = B^2.
\end{eqnarray}
After minor calculation we have
\begin{equation}\label{eta2}
\dot{\gamma}^{p,0} =  \frac{\eta z_0}{\sqrt{z_0^2 - z^2}},\quad \mbox{where}\quad
\eta = 2\sqrt{\bigl(u_0'\bigr)^2 +
\bigl(\frac{u_0}{r}\bigr)^2 + u_0'\frac{u_0}{r}}.
\end{equation}
Substituting \eqref{eta2} into \eqref{gam0}, we obtain the equation for $\frac{\partial u^{p,1}}{\partial z}$.
Solving this equation, integrating with respect to $z$ and matching the first order velocities $u^{s,1}$ \eqref{sh_1ord} and $u^{p,1}$ at $z = z_0$, we get
\begin{equation}\label{u_1p}
u^{p,1} = \eta \sqrt{z_0^2 - z^2}-p_1'(r) (1-z_0)^{\frac{1}{n}}(1+n z_0)- 
g(r)(1-z_0)^{\frac{1}{n}}
.
\end{equation}
Inserting \eqref{mur_pl1}, \eqref{eta2} into \eqref{tau} we find that
\begin{equation}\label{sigmar1}
\tau_{rr}^{p,-1} = \frac{2B }{\eta z_0}u_0'\sqrt{z_0^2 - z^2}, \quad
\tau_{\theta\theta}^{p,-1} =  \frac{2B }{\eta z_0} \frac{u_0}{r}\sqrt{z_0^2 - z^2}.
\end{equation}
Substituting \eqref{sigmar1} into \eqref{pl2}, integrating the resulting equation 
and enforcing continuity of the pressure $p^{s,1}$ \eqref{sh_1ord} and  $p^{p,1}$ at the pseudo-yield surfaces gives
\begin{equation}\label{p100}
p^{p,1} = -\frac{2B}{\eta z_0}\Bigl(u_0' + \frac{u_0}{r}\Bigr)\sqrt{z_0^2 - z^2} + p_1(r).
\end{equation}
We substitute \eqref{sigmar1}, \eqref{p100} into 
\eqref{pl1}, integrate with respect to $z$, using
$\tau_{rz}^{p,1}\left.\right|_{z=0} = 0$, and find
\begin{eqnarray}\label{pl6}
\tau_{rz}^{p,1} = - B\Bigl(z\sqrt{z_0^2 - z^2} + z_0^2 \arcsin\frac{z}{z_0}\Bigr)
\Bigl[\frac{d}{dr}  \Bigl(\frac{2 u_0' + \frac{u_0}{r}}{\eta z_0} \Bigr) +
\frac{(u_0' - \frac{u_0}{r})}{r \eta z_0} \Bigr] \notag\\
 - 2 B \Bigl(\frac{2 u_0'
 + \frac{u_0}{r}}{\eta z_0}\Bigr)z_0'z_0\arcsin\frac{z}{z_0} + p_1' z.
\end{eqnarray}
Enforcing continuity of the  first order shear stress $\tau_{rz}^{s,1}$ \eqref{sh_1ord} and 
$\tau_{rz}^{p,1}$ \eqref{pl6} at the pseudo-yield surface leads to
\begin{equation}\label{match_tau3}
g(r)= - \frac{\pi B}{2}
\Bigl[\frac{d}{dr} \Bigl(\frac{z_0}{\eta} \bigl( 2 u_0' +
\frac{u_0}{r}\bigr) \Bigr) +
\frac{z_0}{\eta r}\Bigl(  u_0'- \frac{u_0}{r} \Bigr)\Bigr].
\end{equation}
To find $p_1$, we insert \eqref{u_1s}, \eqref{u_1p} into the flow rate constraint \eqref{q}, integrate and find
\begin{equation}\label{p21b}
\frac{z_0^2 \eta \pi}{4}-p_1'\Bigl[\frac{B}{z_0}\Bigr]^{\frac{1}{n}-1}\frac{ (2n z_0^2+2n z_0 +n+1)}{(n+1)(2n+1)}-g \Bigl[\frac{B}{z_0}\Bigr]^{\frac{1}{n}-1}\frac{(n z_0+1) }{(n+1)}+u_b^1=0.
\end{equation}
Using \eqref{p21b}, \eqref{dz0ns} we get the expression for $p_1'$:
\begin{equation}\label{p1ns}
p_1' =- \frac{\pi B \eta z_0'}{2} - 2 g u_0'.
\end{equation} 
\begin{figure}[htb]
\includegraphics[height=0.4\hsize,width=0.8\hsize]{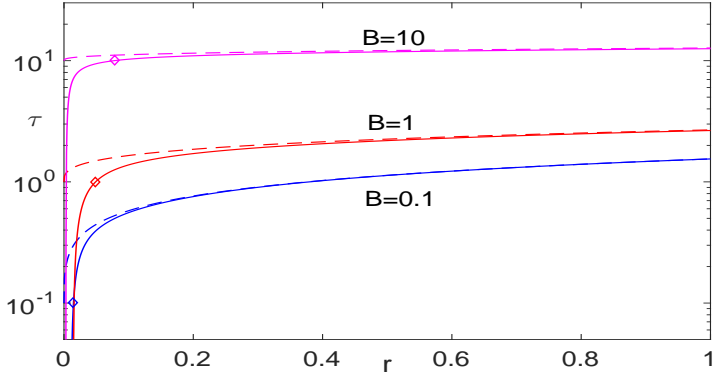}
\caption{
The second invariant of the stress tensor on the upper plate for $\varepsilon=0.1$ and $n=0.5$: dashed lines $\tau^0 =B/z_0(r)$, solid lines $\tau^0 + \varepsilon\tau^1 $ \eqref{tau_1} . } 
\label{mur_fig3}
\end{figure}
Smyrnaios and Tsamopoulos \cite{ST}
showed that unyielded material could only exist around the
two stagnation points of flow. 
The numerical modeling for Bingham   fluid \cite{ST}, \cite{MatMit}, \cite{Mur_pof}, experiments \cite{Lanos06} and asymptotic solution \cite{Mur_JNNFM2017} confirmed the results \cite{ST}.
It is interesting to investigate the obtained asymptotic solution for Herschel-Bulkley fluid near the stagnation point. We will examine the second invariant of the stress $\tau$ at the plate at $z = 1$.
\begin{equation}
\label{tau_1}
\tau(r,1) = |\tau^0_{rz}(r,1) + \varepsilon \tau^1_{rz}(r,1)| + O(\varepsilon^2)=\frac{B }{z_0(r)} - \varepsilon \Bigl(p_1' (r) + g(r)\Bigr) + O(\varepsilon^2 ),
\end{equation}
where $p_1'(r)$ and $g(r)$ are determined by \eqref{match_tau3},
\eqref{p1ns}. 
The leading order $\tau^0=\frac{B }{z_0(r)}$  
exceeds the yield stress, because of $z_0(r) \le 1$,
hence, according to the zero-order solution, the Herschel-Bulkley fluid in the shear region is yielded.
To analyze the behavior of the function $\tau_{rz}(r,1)$ near the axis of symmetry, we  expand function $z_0(r)$ in Taylor series.
The expressions for $\tau_{rz}^{0}(r,1)$ \eqref{mur_tau_s} and $\tau_{rz}^1(r,1)$ \eqref{sh_1ord}, \eqref{match_tau3}, \eqref{p1ns} takes the form:
\begin{equation*}
\tau^0_{rz}(r,1)=-B-\Bigl[\frac{r B(n+1)}{2n }\Bigr]^{\frac{n}{n+1}}- \Bigl[\frac{B(n+1)}{n}\Bigr]^{\frac{n-1}{n+1}}\frac{(3n+1)}{(2n+1)}\Bigl[\frac{r}{2 }\Big]^{\frac{2n}{n+1}}+\mathcal{O}({r^{\frac{3n}{n+1}}}),
\end{equation*}
\begin{equation}\label{tau_ser}
\tau^1_{rz}(r,1)=\frac{\sqrt{3}\pi B^{n/(n+1)}}{4}
\Bigl[\frac{2 n}{r(n+1)}\Bigr]^{1/(n+1)} +\mathcal{O}\Bigl({r^{(n-1)/(n+1)}}\Bigr).
\end{equation}
The expression in the brackets on the right-hand side of \eqref{tau_ser} tends to infinity when $r \rightarrow 0$. In the limit $r\rightarrow 0$, we have, 
$\tau_{rz}(r,1) = -\frac{B }{z_0(r)} + \varepsilon \Bigl(p_1' (r) + g(r) \Bigr) \rightarrow  \infty$.
Therefore, about $r = 0$ there is a point $r = r_0$, at which $\tau_{rz} = - B$ and $\tau = B$.
For $r \le r_0$ we have $\tau < B$.
So, the asymptotic analysis predicts a region of unyielded material near the axis of symmetry.
Consequently, our asymptotic expansion breaks down.

Fig.~\ref{mur_fig3}
shows that the values of the second invariant $\tau$ of the stress tensor along the disk surface increase with the Bingham number $B$. The solid lines, corresponding to
  $\tau(r,1) = \tau^0(r,1) + \varepsilon \tau^1(r,1)$, is below the dashed line, corresponding to $\tau^0(r,1)$. The function $\tau(r)$ 
 increases with $r$ and decreases sharply near the axis of symmetry $r =0$. The points on the graphs for which $\tau = B$ are marked diamonds on the solid lines. 

\subsection{Pressure distribution and squeeze force}\label{sec_force}
In order to obtain the squeeze force we first calculate the normal radial stress $\sigma_{rr}$, neglecting the terms $O( \varepsilon^2)$:
\begin{equation}\label{sigmar3}
\sigma_{rr}(r,z) =
\begin{cases}
 - p_0(r)-\varepsilon p_1(r),\quad &  z \in (z_0,\,1],\\
 - p_0(r)-\varepsilon p_1(r) +\frac{2B }{\eta z_0}(2u_0'+ \frac{u_0}{r})\sqrt{z_0^2 - z^2}\quad, &z \in [0,\,z_0].
\end{cases}
\end{equation} 
 We see that $\sigma_{rr}$ is dependent on $z$, which makes it impossible
 to satisfy the zero-stress condition  at $r = 1$ exactly. We impose the average boundary condition, suggested in \cite{SherDur96}:
\begin{equation}\label{str_av}
 \int_{0}^{1} {\sigma}_{rr}(1,z)~dz = 0.
\end{equation}
We  substitute \eqref{sigmar3}
into the average condition \eqref{str_av} and obtain
\begin{equation}\label{pconst}
p_0(1) + \varepsilon p_1(1) = \varepsilon p_R = \varepsilon B \pi \Bigl[\frac{ z_0(1)(2 u_0'(1) + u_0(1))}{2\eta(1)}\Bigr].
\end{equation}

The pressure gradient zero-order is given by  $p_0'(r) = - \frac{B}{z_0}$.
The integration, using the results from 
\eqref{buck}, \eqref{dz0ns}:
\begin{equation}\label{r_dr1}
r = \frac{2 B^m(1-z_0)^{m+1}(m +1+z_0)}{z_0^m (m+1)(m+2)},
\end{equation}
\begin{equation}\label{r_dr2}
dr = -\frac{2 B^m(1-z_0)^m(2z_0^2 +2m z_0+m+m^2)}{z_0^{m+1}(m+1)(m+2)}dz_0,
\end{equation}
produces the following pressure distribution ($p_0(1) = p_R$):
\begin{equation}\label{p_s}
p(r)=2 B^{m+1}\int_r^1 
\frac{(1-z_0)^m(2z_0^2 +2m z_0+m+m^2)}{z_0^{m+2}(m+1)3(m+2)}dz_0 +\varepsilon p_R.
\end{equation}
The function $p_1(r)$ can
be found from \eqref{p1ns} numerically with $p_1(1)=0$.

\begin{figure}[htb]
(a)\includegraphics[height=0.36\hsize,width=0.48\hsize]{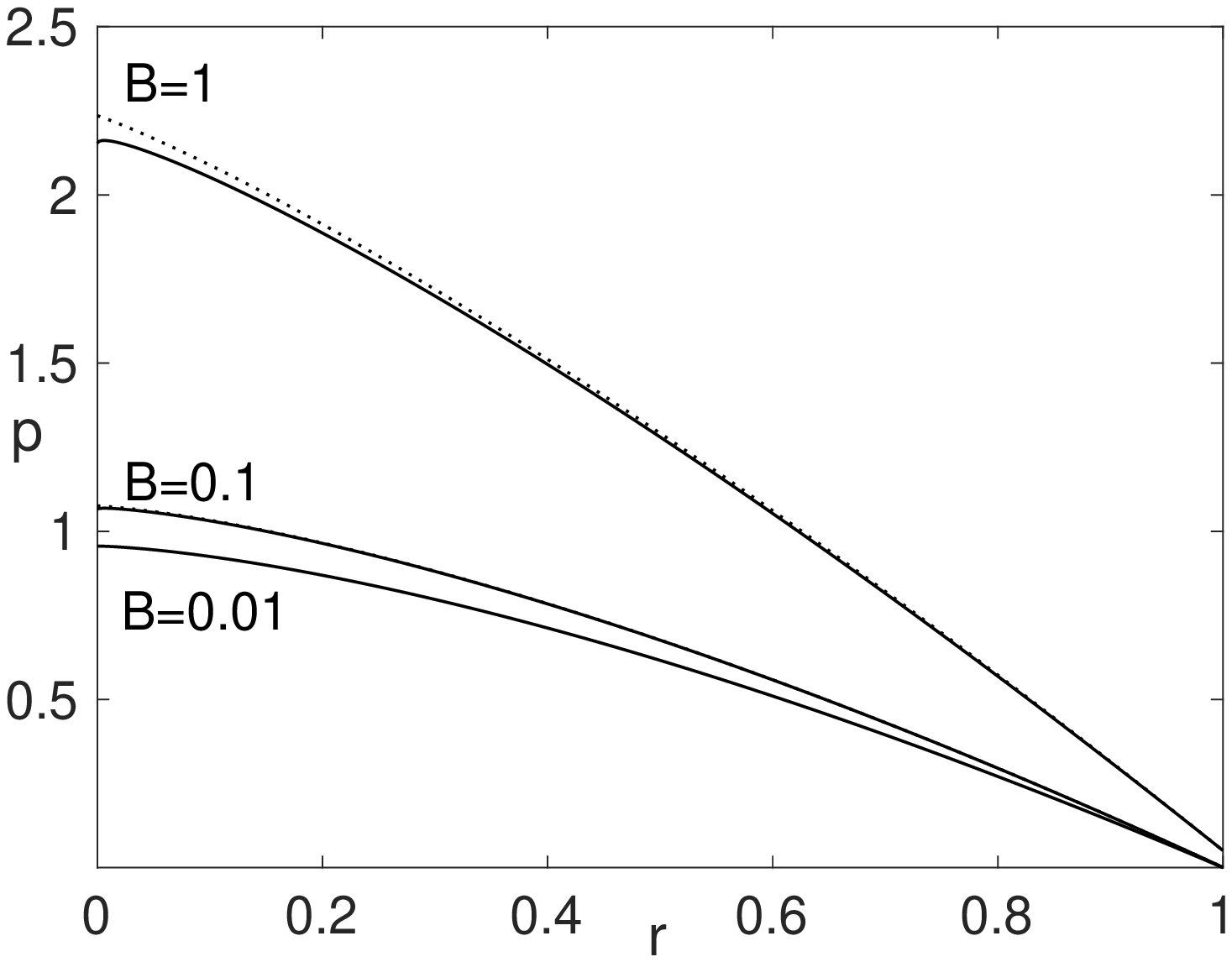}
(b)\includegraphics[height=0.36\hsize,width=0.48\hsize]{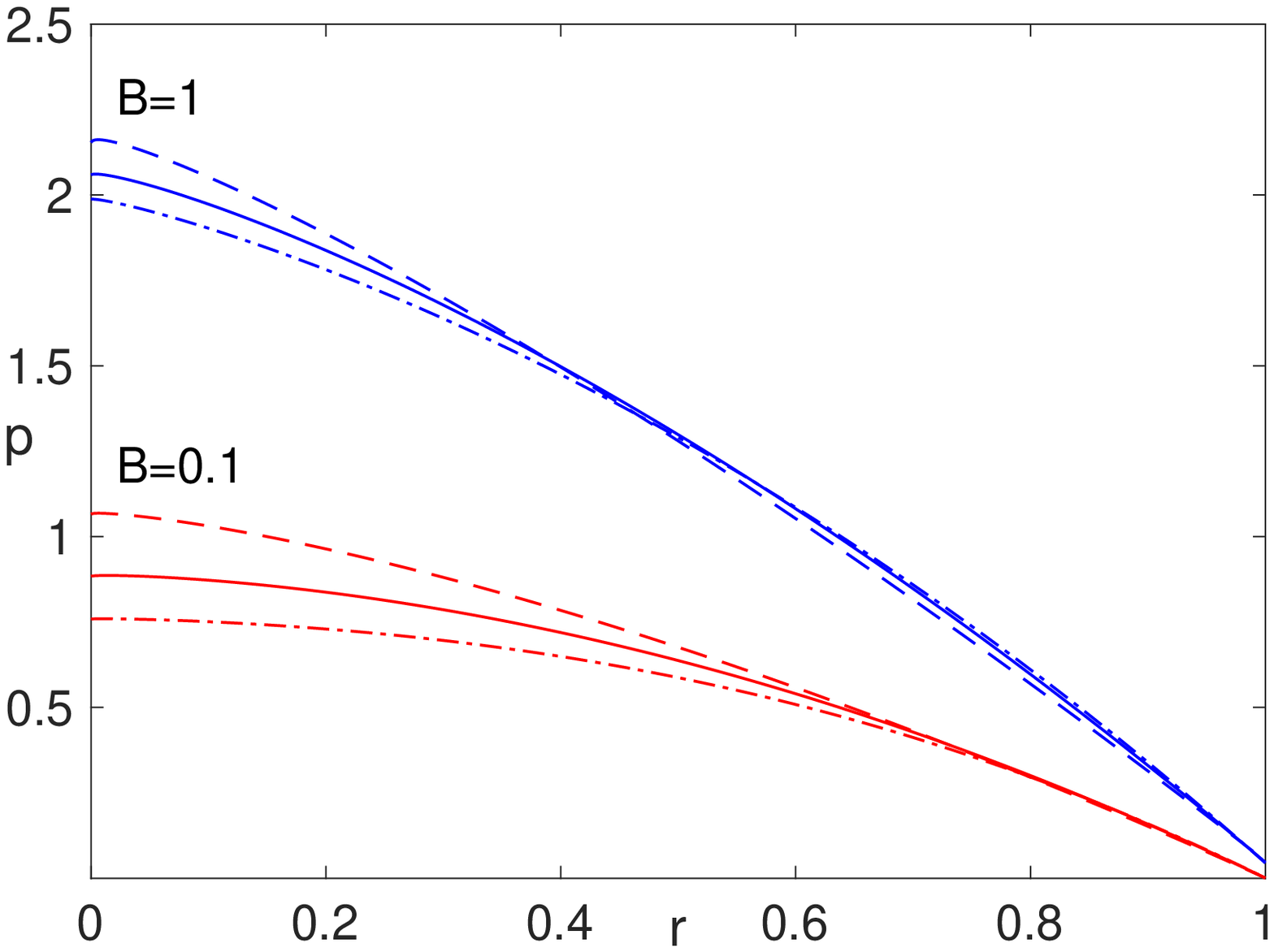}
\caption{The pressure distribution on the upper
plate for $\varepsilon = 0.1$: (a) for $n=0.5$, the solid lines indicate the pressure $p^{s}=p_0+\varepsilon p_1$, \eqref{pconst}, \eqref{p1ns}, \eqref{p_s}, the dotted lines indicate the pressure zeroth order $p^{s,0}=p_0$ \eqref{pconst}, \eqref{p_s}; (b) for $n=0.5$ (dashed lines), $n=1$ (solid lines), $n=1.5$ (dash-dotted lines).}
\label{mur_fig4}
\end{figure}

The solid lines show computed final states, the dotted lines denote the leading-order result (8) and the dashed line shows the higher-order prediction (9) .

Fig.~\ref{mur_fig4}(a)
shows the pressure distribution along the disk surface for different
values of $B$ .
The values of $\varepsilon p_R$ are very small for $B = 0.01, 0.1$,
increase with the Bingham number. 
We see that adding $\varepsilon p^{s,1}$ reduces the value of pressure, so the solid line ($p^{s}$) is lower than the dotted line ($p^{s,0}$).
Since for $B = 0.01$, $0.1$ the value of $\varepsilon p^{s,1}$ is very small additive to the $p^{s,0}$,
the plots of $p^{s}$ and $p^{s,0}$ coincide. We see in Fig.~\ref{mur_fig4}(b) that as power-law index $n$ decreases, the pressure becomes more significant in the center of the disk and smaller at the edge of the disc. With the increase in the number of Bingham, the graphs for different power-law indices $n$ less differ among themselves.

\begin{figure}[htb]
\includegraphics[height=0.45\hsize,width=0.8\hsize]{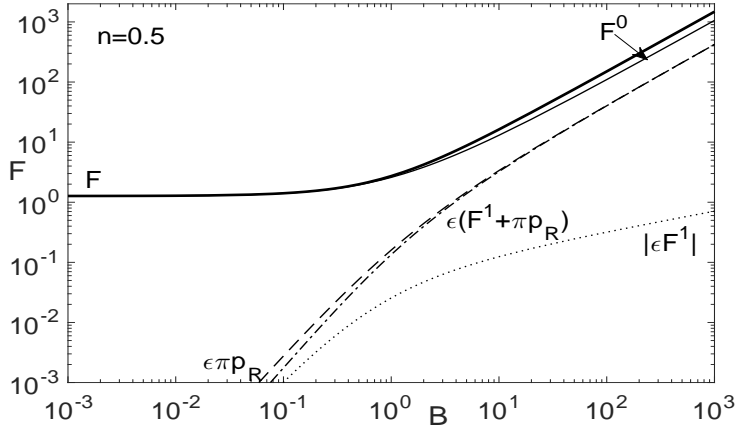}
\caption{
The squeeze force as a functions of Bingham number $B$ for $\varepsilon=0.1$:
(a) $F=F^0+\varepsilon(F^0+p_R)$, $F^0$ \eqref{F_0}, $|\varepsilon F^1|$ \eqref{F_1}, $\varepsilon \pi p_R$ \eqref{pconst}.}
\label{mur_fig5}
\end{figure}
We calculate the squeeze force by integrating
the normal axial stress $\sigma_{zz}^s = - p^s(r)+ \varepsilon^2 \tau_{zz}^{s,0}(r,z)$ over the plate surface. In the lubrication solution up to second order in $\varepsilon$ the component $\tau_{zz}$ is negligible. So
\begin{eqnarray}\label{Fn0sl}
&F =2 \pi \int_0^1 p^{s} ~ r~dr =
 ( \pi r^2 p^{s} )\left.\right|_0^1 - \int_0^1 \pi r^2 \frac{dp^{s}}{dr}dr
= F^0 + \varepsilon F^1 + \varepsilon \pi p_R,
\\\label{F0_F1}
&F^0 = -\pi \int_0^1  p_0' r^2~dr,\quad
F^1 = -\pi \int_0^1 p_1' r^2~dr.
\end{eqnarray}
We calculate $F^0$, using \eqref{mur_tau_s}, \eqref{F0_F1}, \eqref{r_dr1}, \eqref{r_dr2}, 
\begin{equation}\label{F_0}
F^0 = 
8 \pi B^{3m+1}\int_0^1 
\frac{(1-z_0)^{3m+2}(m +1+z_0)^2(2z_0^2 +2m z_0+m+m^2)}{z_0^{3m+2}(m+1)^3(m+2)^3}dz_0.
\end{equation}
To calculate $F^1$, we
substitute the expression for $p_1'$
\eqref{p1ns} into \eqref{F0_F1},
integrate by parts , and obtain
\begin{equation}\label{F_1}
F^1 = \pi^2 B \Bigl[ \frac{  z_0(1) u_0(1)}{\eta(1)} (2 u_0(1) + u_0'(1)) - \frac12 \int_0^1 \eta z_0 r~dr \Bigr].
\end{equation}
Fig.~\ref{mur_fig4}
suggests that there is very little difference between the zeroth order and first order results for the pressure over the plate. As a result, there will be very little difference between the zeroth and first order forces.
In Fig.~\ref{mur_fig5}
we plot separate results for the total force $F$, $F^0$, $|\varepsilon F^1|$, $\varepsilon p_R$ as the functions of $B$. (Since $F^1$ is negative and we use a logarithmic scale, we plot $|\varepsilon F^1|$.) All these quantities increase with increasing B, the function $\varepsilon p_R(B)$ grows much faster than $|\varepsilon F^1(B)|$. For $B < 1$ the graphs $F$ and $F^0$ coincide, for $B > 1$ the graph of $F$ is higher than the graph of $F_0$.
The total force $F$ is larger than $F^0$ due to $\varepsilon p_R$, 
and this difference increases with $B$. 

\begin{figure}[htb]
(a)\includegraphics[height=0.35\hsize,width=0.8\hsize]{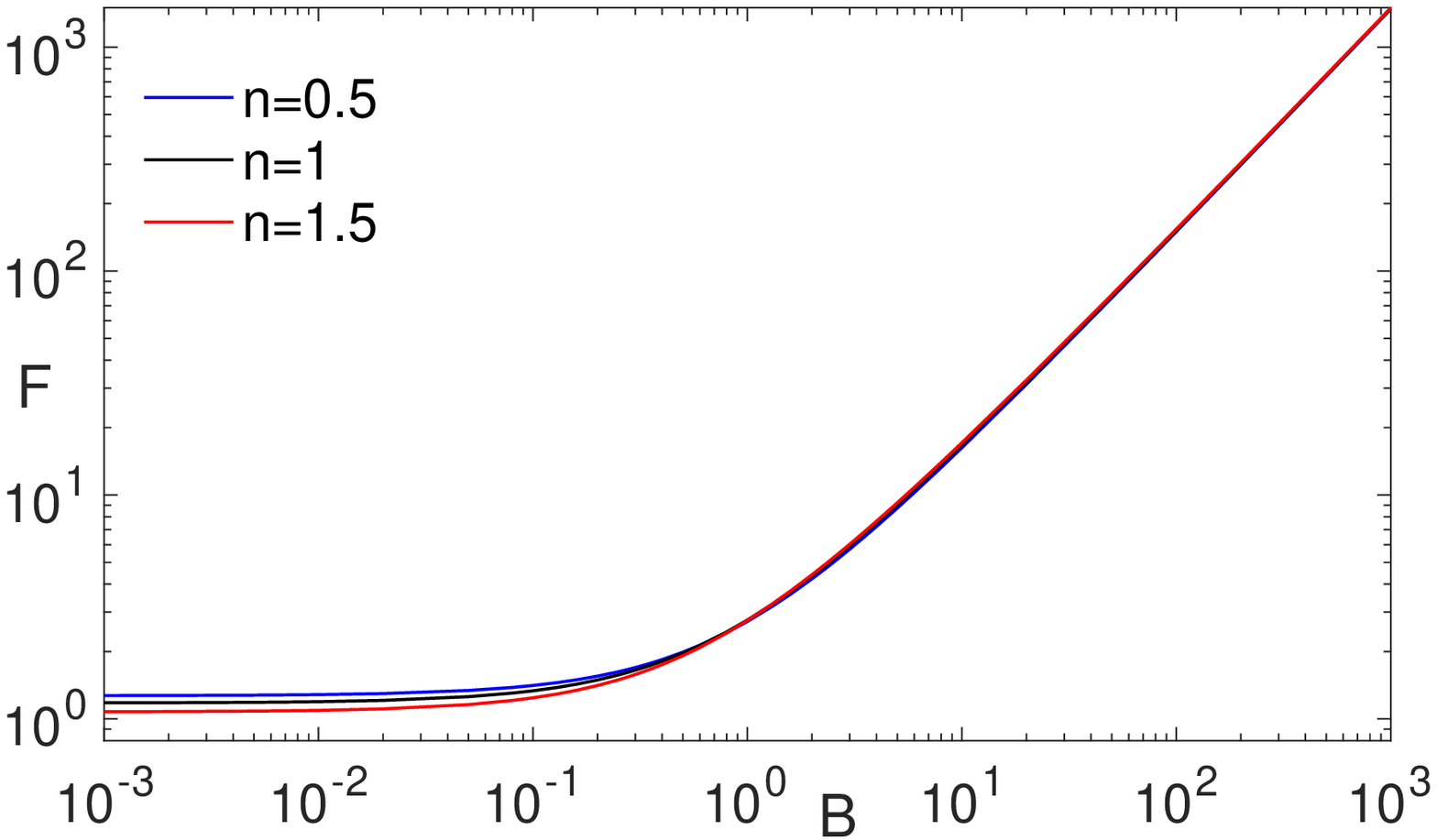}
(b)\includegraphics[height=0.35\hsize,width=0.8\hsize]{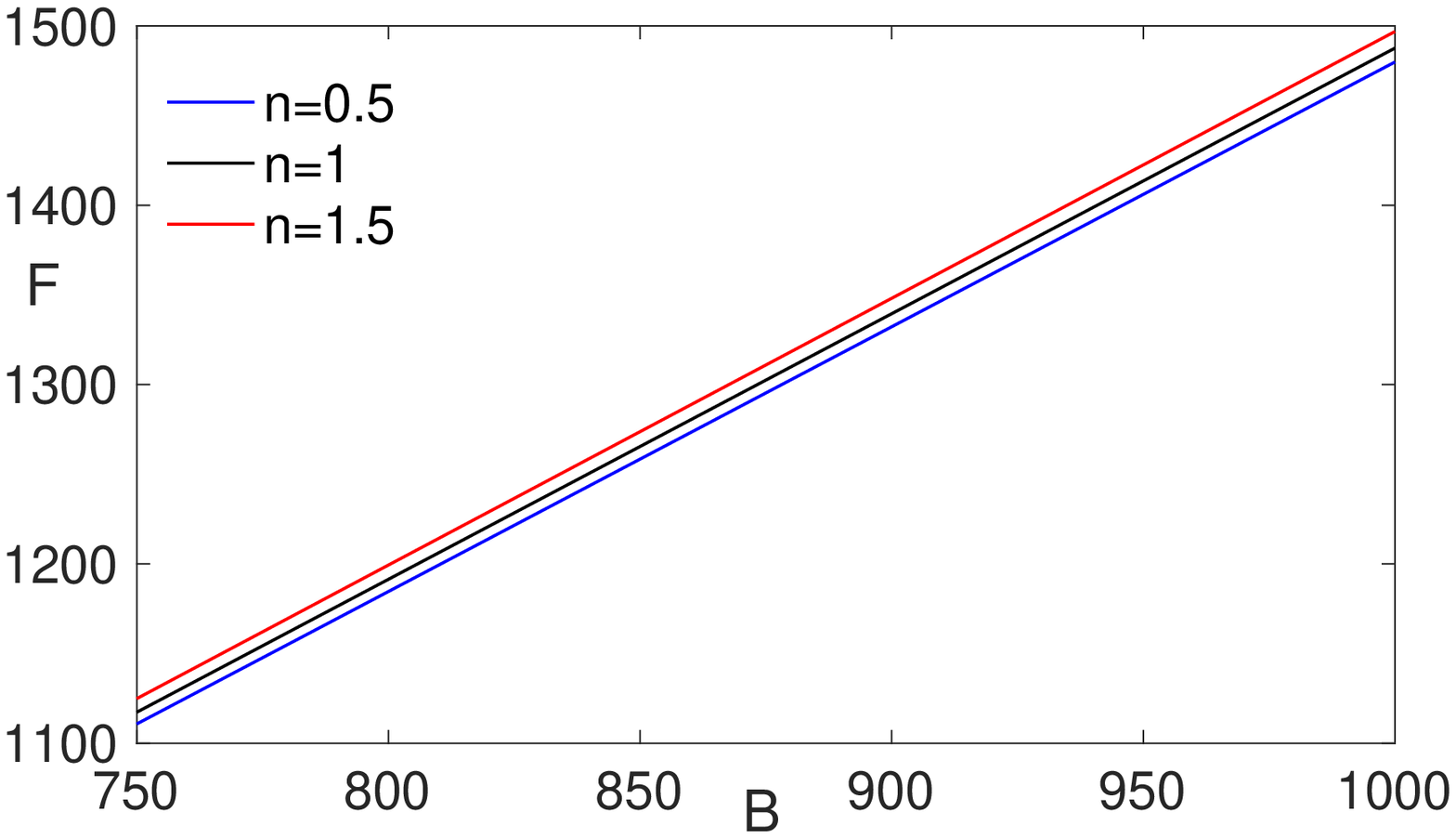}
\caption{
The squeeze force as a functions of Bingham number $B$ and power index $n$ for $\varepsilon=0.1$: 
(a) a logarithmic scale,
(b) a linear scale.
}
\label{mur_fig6}
\end{figure}
Fig.~\ref{mur_fig6} shows the dependence of the total force on the Bingham number and the power index. We observe an interesting effect: for small Bingham numbers, the total force greater for low power indices; however, but with increasing Bingham number, the total force becomes greater for higher power indices. Because the graphs in a logarithmic scale are the same for large values of B, in Figure 6(b), we show the plots for large Bingham number in a linear scale.

\end{document}